\documentclass{amsart}[12pt]
\usepackage{amsmath,amsfonts,amssymb,latexsym,amsthm}
\usepackage{graphicx,epsf}
\usepackage{epsfig}
\usepackage{psfrag}

\def\zz{\mathbb Z}

\def\cc{\mathbb C}

\def\Ga{{{\rm G}}}

\def\De{\Delta}

\def\al{\alpha}

\def\wt{\widetilde}
\def\<{\langle}
\def\>{\rangle}

\def\ts{\hskip.015cm}

\def\0{{\mathbf 0}}

\def\sS{{\tt S}}

\def\rA{z}

\begin{document}
\title{The area of cyclic polygons: \, Recent progress on
Robbins' Conjectures}

\author[Igor~Pak]{ \ Igor~Pak$^\star$}


\thanks{\thinspace ${\hspace{-.45ex}}^\star$Department of Mathematics,
MIT, Cambridge, MA, 02139.
\hskip.06cm
Email:
\hskip.06cm
\texttt{pak@math.mit.edu}}
\maketitle

 {\hskip5.83cm
 July 30, 2004
 }

\vskip1.cm

\begin{abstract}
In his works~\cite{R1,R2} David Robbins proposed several
interrelated conjectures on the area of the polygons inscribed
in a circle as an algebraic function of its sides.  Most recently,
these conjectures have been established in the course of
several independent investigations.  In this note we give an
informal outline of these developments.
\end{abstract}

%

\vskip.9cm

\section{Introduction}\label{intro}

Let $a_1, a_2, \ldots, a_n$ be the side lengths of a convex
polygon inscribed in a circle.  What is the area
$\sS = \sS(a_1,a_2,\ldots,a_n)$ of the polygon as
a function of the sides?  This question goes back
to Heron of Alexandria (the case $n=3$) and
Brahmagupta (the case $n=4$).  It seems, David
Robbins was the first to address this question in
full generality and suggest the way of phrasing the
answer~\cite{R1,R2}.  We start with the general remarks
on the problem (largely due to Robbins)
and then outline recent developments in an informal
essay style.

\smallskip

Following Robbins, we call polygons inscribed in a
circle the {\em cyclic polygons}. We denote the vertices
by $A_1,A_2,\ldots,A_n$ and the center by~$O$.

First, observe that $\sS(a_1,a_2,\ldots,a_n)$ is well defined,
that is there exists at most one cyclic polygon with the given
(ordered list of) side lengths.
Indeed, start with a large enough circle and place $n+1$
vertices $A_0,A_1, \ldots,A_n$ at distance $|A_i-A_{i+1}| = a_i$.
Continuously decreasing the radius we obtain a unique
convex polygon with $A_0 = A_n$, as desired.

Second, observe that $\sS$ is a symmetric function in~$a_i$.
Indeed, this follows from the fact that we can interchange
triangles $[OA_{i-1}A_{i}]$ and $[OA_iA_{i+1}]$.  The details
are straightforward.

Our final observation is the fact that $\sS$ is an algebraic
function of the side lengths~$a_i$.  First, notice that it is
polynomial in the coordinates of vertices~$A_i=(x_i,y_i)$:
$$(\circ) \hskip1.cm
\sS \, = \, \frac12 \, \left|
\begin{matrix}
x_1 & x_2 \\
y_1 & y_2
\end{matrix}\right|
\, + \, \frac12 \,
\left|
\begin{matrix}
x_2 & x_3 \\
y_2 & y_3
\end{matrix}\right|
\, + \ \ldots
\, + \, \frac12 \,
\left|
\begin{matrix}
x_n & x_1 \\
y_n & y_1
\end{matrix}\right|. $$
Here each summand is equal to the (signed) area of the
triangle $[\0 A_i A_{i+1}]$, and it is easy to see that
they add up to the area of the polygon ($\0$ denotes the origin).

Now, move the polygon in such way that~$A_1 = (0,0)$
and $A_2 = (a_1,0)$. There are $2(n-2)$ free variables for
the remaining vertex coordinates and~$2$ variables for the
coordinates of the center~$O$. Together these give~$2n-2$
variables. Similarly, the remaining side lengths give~$n-1$
equations, the equality of the distance to the origin give
another~$n-1$ equations, which total~$2n-2$ equations.
One can show that these equations are algebraically
independent so all free coordinates are in fact the
algebraic functions of the side lengths~$a_i$.  Thus,
from the formula~$(\circ)$, so is the area~$\sS$.

Note that depending on the orientation of the polygon,
the (signed) area~$\sS$ given by~$(\circ)$ is either
positive or negative. Also, in the quadratic equations
for the distances
$$(x_{i} - x_{i+1})^2 + (y_{i} - y_{i+1})^2 = a_i^2
$$
only the squared edge lengths appear.
Thus, for each~$n$ there exist a minimal polynomial equation
$P_n(\sS^2,a_1^2,\ldots,a_n^2)=0$ which has the squared
area as its root.  Changing the first variable, we obtain
$\al_n(16\ts\sS^2,a_1^2,\ldots,a_n^2) = 0$,
which Robbins called the {\em generalized Heron polynomials}.

As we mentioned above, polynomials $\al_3,\al_4$ were well known.
In his work~\cite{R1,R2}, Robbins calculated $\al_5, \al_6$ and
made a number of conjectures on the general form of
polynomials~$\al_n$.  By now, his conjectures have largely been
established in a series of recent development.  Before we move
on to outline their solutions, let us mention that Varfolomeev,
unaware of Robbins' work, recently rediscovered some of his
results and independently made a number of advances on the
subject~\cite{V1,V2}.

As the reader will see, we do not include any technical
details, nor do we present a formal survey.  Instead,
give the reader a quick introduction to the subject
and its basic ideas, aiming to ease the entrance barrier
and to simplify navigation through recent developments.

\vskip.9cm

\section{The first coefficient}\label{first}

Based on small examples, David Robbins conjectured that
the polynomials $\al_n$ are {\em monic} in
variable~$\rA = 16 \ts\sS^2$ (``monic'' means that the
highest coefficient is equal to one).
This seemingly random observation is in fact very interesting
and is strongly related to Sabitov's theory and the proof of
the bellows conjecture.  This connection was independently
discovered by Connelly~\cite{C3} and Varfolomeev~\cite{V1},
who gave two different proofs.  Let us first elaborate on
Sabitov's work.

The story goes back to Connelly's celebrated discovery of
flexible (nonconvex) polyhedra and his {\em bellows
conjecture}\footnote{
Bob Connelly declines to take credit for the bellows
conjecture and wrote to me that it was communicated to
him by various people, all of whom refer yet again to
other people.  Therefore, the conjecture is a folklore
in the area, while Connelly deserves a great deal of credit
for its advancement.
}
as to whether the volume remains invariant under flexing
(continuous face-preserving deformations).  We refer to~\cite{C1}
for background and references.  Later both Connelly
and Sabitov conjectured that in fact the volume is integral
over the ring generated by squares of polyhedra's edge lengths.
This immediately implies the conjecture since nonzero polynomials
have only finitely many roots, and thus allow only a finite
set of possible volume values.

The bellows conjecture was established by Sabitov, who gave
several consequently improving expositions of his proof in
a series of papers (see e.g.~\cite{S1,S2}).  Let~$\ell_i$
denote the edge lengths of the polytope with edge graph~$\Ga$,
and let~$V$ denote its volume.  Sabitov showed that
there exists a nontrivial polynomial equation
$P_\Ga(V^2,\ell_1^2,\ell_2^2,\ldots)=0$.
The difficult part in the proof is not computing
this polynomial but checking that the leading coefficient
is not zero.  In fact, after a change of variables
$\wt P_\Ga(144\ts V^2,\ell_1^2,\ell_2^2,\ldots)=0$ all
coefficients become integral, and the polynomial
$\wt P_\Ga$ is monic in~$(144\ts V^2)$.

Unfortunately, Sabitov's proof is based on elimination
theory and is more technical than enlightening.
Sabitov's approach was later modified in~\cite{CSW}
where the \emph{theory of places} is used to prove the
bellows conjecture.

Note that when the polytope is a simplex the resulting
polynomial equation can be viewed as (a different) generalization
of the Heron formula~\cite{S1}.  The striking similarity
of two problems led Varfolomeev to rediscovery of some of
Robbins' ideas and results.  He used Sabitov's methods to
show that polynomials~$\al_n$ are monic~\cite{V1}.

Connelly came to his proof~\cite{C3} independently,
after~\cite{wsj}
advertised Robbins' efforts.  He  observed the similarity
as well, and used the theory of places, to obtain a
beautiful and concise proof of this Robbins conjecture.

In the spirit of the bellows conjecture, both authors
address the question as to when cyclic polytopes are
flexible.  The immediate implication of the above result
is the fact that the (symplectic) area is unchanged under flexing.
In fact, as was shown by Connelly much earlier~\cite{C1},
the  area of flexible cyclic polygons is always zero.
This was also rediscovered by Varfolomeev~\cite{V1}
(see also~\cite{C3}).

\vskip.9cm

\section{The degree}\label{degree}

The most aesthetically pleasing conjecture of Robbins is his
proposed formula for the degree of generalized Heron polynomials
(in the variable $\rA=16\ts\sS^2$):
$$(\divideontimes) \ \ \ \,
\deg \ts \al_{2k+1}  = \De_k, \ \, \deg \ts \al_{2k+2}  = 2 \De_k, \ \
\text{where} \ \De_k = \, \frac{2k+1}{2}\binom{2k}{k} \, - \, 2^{2k-1}.
$$
Let us explain the origin of this formula for~$\De_n$.  Observe
that formula $(\circ)$ works not only for (the usual)
cyclic polygons, but also for those with self-intersections,
still inscribed in a circle.  Therefore, the minimal polynomial
$\al_n (\rA)$ has at least as many real roots as the number of
different areas of these self-intersecting polygons.  Then,
Robbins showed that for nearly equal side lengths $a_1,a_2,\ldots,a_n$
all self-intersecting cyclic polygons have different areas.
A simple combinatorial argument gives the r.h.s.
in~$(\divideontimes)$ for the number of different
self-intersecting cyclic polygons, and implies
the desired lower bound on the degree of~$\al_n$.

Robbins' conjecture for the degree, the formula~$(\divideontimes)$,
was established in~\cite{FP}, and later by a simpler but related
argument in~\cite{MRR}.  Both proofs first obtain the corresponding
formula for the degree of polynomial equations on the radius of
the circle, and then move to the degree of~$\al_n$.
The study in~\cite{FP} goes much deeper, as the authors establish
formal connections with Sabitov's theory, which we outline below.

It was observed by Sabitov that not only the volume, but other
``polynomial invariants'' of polytopes are roots of polynomial
relations with coefficients being polynomial in the squared side
lengths.  It is a natural general question to compute the degree
of these minimal polynomials.  We should emphasize that here we
discuss only convex polytopes, so certain technical difficulties
of Sabitov's approach do not appear in this case.  Now, for convex
polytopes not only the volume, but {\em all} diagonals are the roots
of polynomial equations and thus one can think of them as of an
extension of the Cauchy rigidity theorem.  Following this logic,
in~\cite{FP} we refined a known algebraic proof of the Cauchy rigidity
theorem and added an argument from algebraic geometry. We obtained
a general upper bound on the degree of minimal polynomial
relations for all polynomial invariants of convex polytopes,
including the volume, (squared) diagonal lengths, etc.   In this
special case the upper bound we obtain for $P_\Ga$ is in terms
of {\em complex} realizations of the graph~$\Ga$
(realizations in~$\cc^3$)  and
in the worst case gives~$2^m$, where~$m$ is the number of edges
in~$\Ga$ (= the number of edges in the polytope).

Now, after we learned from~\cite{wsj} about Robbins' conjectures,
we discovered a formal connection between our work (Sabitov's
theory) and that of Robbins.  Consider a bipyramid with
(large enough) equal
length edges leaving north or south poles, and the edge lengths
$a_1,a_2,\ldots,a_n$ in the middle.  Clearly, the middle edges
form the desired cyclic polygon, and in fact different (real)
realizations of this bipyramid produce different
(self-intersecting) cyclic polygons.  Also, the {\em main} (north to
south pole) {\em diagonal} is related by Pythagoras theorem to the
radius.  It may seem like our main upper bound is directly
applicable in this case to obtain the degree of the minimal
polynomial relation for the radius.  The problem is that the
number of complex realizations is a difficult quantity to
compute in most cases.  In fact, our logic moves backwards
and is more convoluted.

First, we use an ad hoc combinatorial argument to compute explicitly
the minimal polynomial relation for the radius and its degree.
The corresponding polynomial relation turns out to have a nice
closed formula amenable to direct calculation.  Then we use
the relationship described above to obtain the polynomial
relation for the main diagonal, and thus bound the number of
complex realizations, which we show is equal to the number of
real realizations (self-intersecting cyclic polygons).
Finally, we apply our upper bound
theorem to obtain the upper bound on the degree of~$\al_n$,
the minimal polynomial for the area of cyclic polygons.  With
Robbins' matching lower bound we obtain the result.

It is interesting to note that we never actually obtain any
useful formula for~$\al_n$.  In fact, the only polynomial
invariant that stands out in this case is the radius---all
others play a supporting role.  For example, instead of the
area we could be proving formula~$(\divideontimes)$ for the
degree of the minimal polynomial relation on the sum or
squares of all diagonals in a cyclic polygon.

Let us say a few words on the proof in~\cite[$\S 5$]{MRR}.
This work started out by David Robbins and Julie Roskies
just a few month before Robbins' premature death, and continued
later with the help of Miller Maley.  Their proof of the degree
formula
starts with the use of M\" obius formulas for the radius of the
cyclic polygon, which are essentially equivalent to our formulas.
Rather than utilize the general theory we develop in~\cite{FP},
the authors use a simplified algebraic geometry argument adjusted
in this particular case to obtain the result.  Basically,
their argument is the same as our argument for the special
case of a bipyramid.

Finally, let us note that Varfolomeev also studies explicit
formulas for the radius in cyclic polygons~\cite{V1}.
He also guesses the answer in terms of self-intersecting
polygons, but never obtains a general formula nor even
calculates their number beyond few small cases.

\vskip.9cm

\section{Explicit formulas}\label{closed}

One more Robbins' conjecture concerns the form of polynomials
of~$\al_{2k-1}$ versus that of~$\al_{2k}$.  Roughly speaking,
he claimed that given the formula for~$\al_{2k-1}$ one can easily
obtain the formula for~$\al_{2k}$ as a product of  the formula
for~$\al_{2k-1}$ and its variation.  This conjecture was established
by Varfolomeev~\cite{V1} by a direct argument (see also~\cite{MRR}).
As a corollary, calculations of Robbins et al. for cyclic pentagons
and heptagons immediately translate to give the formulas for
cyclic hexagon and octagons~\cite{R1,MRR}.

It was Robbins' wish to obtain a concise formula for~$\al_7$ and
although he did not live to finish the project, such a formula was
recently obtained in~\cite{MRR}.  Of course, Robbins already showed
that some kind of formula exists, but given the large number
of terms one can ask if there is a way to simplify it.
In view of our earlier impression (see above) that a nice formula
may exist only for the radius, we find it amazing that the authors
were able to obtain a concise formula for the area.

Let us mention that already the Robbins' formula for $\al_5$
is very interesting as it expresses the area
as a discriminant of a certain ``mystery cubic''~\cite{R2}.
It remains unclear where this cubic comes from and what is
its role in the grand scheme of things.
For example, Varfolomeev~\cite{V1} does not notice this
formula and uses rather elaborate explicit formulas for~$\al_5$.

Now, in~\cite{MRR} the authors obtain a closed formula for~$\al_7$
in terms of a resultant of two concise, but not generalizable
polynomials.  The resulting formula is nice but again very
mysterious.  It remains open whether this work can be extended
to obtain concise formulas for~$\al_n$, where~$n \ge 9$.

\vskip.9cm

\section{Final remarks and open problems}\label{final}

From the point of view of Sabitov's theory and our paper~\cite{FP},
it would be natural to ask for the minimal degree polynomial
relations for the volume or the diagonal lengths in various families
of convex polytopes.  As noted in~\cite{FP}, even for relatively
small polytopes this problem is computationally intractable and
new ideas are needed even to obtain the exact asymptotic behavior.
At the moment, the precise formulas that we found for regular
bipyramids are clearly beyond reach in most cases.

An interesting twist on cyclic polygons was proposed in~\cite{MRR}
where the authors define what they call ``semicyclic polygons'',
where one side is forced to be a diameter and its length is not
specified.  It seems that much of the work extends to this case with
little difficulty.  We propose to consider an equivalent model
of the centrally symmetric cyclic polygons with given edge lengths.
This version has the advantage of being possible to generalize
to cyclic polygons with $\zz/\ts k\ts\zz$ cyclic symmetry.  It would be
interesting to see if the analysis extends to this case.
In general, one can consider general polytopes with a given
symmetry group.  Developing the corresponding ``equivariant
Sabitov theory'' seems like a fruitful direction.

When it comes to the area and generalized Heron polynomials~$\al_n$,
it is probably too much to ask for a concise general formula.
Still, we remain optimistic of other research venues.
In his latest work~\cite{V2},
Varfolomeev calculated the Galois group of~$\al_5$ and showed that it
is the full group of permutations~$S_7$ (he did this also for the
radius).  There seem to be no immediate implications of this result
except perhaps the impossibility of ``construction'' of
the cyclic pentagon with a ruler and a compass, given the
generic lengths of edges (such construction of a regular
pentagon is well known).
In any case, it would be nice to extend these calculations for
general~$\al_n$.

Finally, further connections to rigidity theory are waiting
to be explored.  We refer to final remarks in~\cite{C3} for
directions and motivation.  Also, an intriguing construction
of a finitely generated infinite-dimensional Lie algebra
was announced in~\cite{V2} and promised to be the subject
of the future investigations.  We are anxious to see
how this theory further develops.

\vskip.9cm

\subsection*{Acknowledgments} \label{ack}
\
I am very grateful to Lynne Butler, Miller Maley, Julie Roskies,
Izhad Sabitov and Vitalij Varfolomeev for informing me of their
results and the interest in my work.  Special thanks to Bob
Connelly for his help, encouragement and for reading the
previous version of this mini survey.  Finally, I would like
to thank my coauthor and advisee Maksym Fedorchuk for his
bravery in choosing to work in an unfamiliar area with
an unfamiliar advisor.

\vskip.9cm


\end{document}